\documentclass[12pt]{amsart}
\usepackage{color,ulem}
\numberwithin{equation}{section}

\def\1#1{\overline{#1}}
\def\2#1{\widetilde{#1}}
\def\3#1{\widehat{#1}}
\def\4#1{\mathbb{#1}}
\def\5#1{\mathfrak{#1}}
\def\6#1{{\mathcal{#1}}}

\def\C{{\4C}}
\def\R{{\4R}}

\def\Z{{\4Z}}

\def\aut{{\sf aut}}
\def\Aut{{\sf Aut}}

\def\Re{{\sf Re}\,}
\def\Im{{\sf Im}\,}

\def\hol{{\sf hol}}
\begin{document}

\def\codim{{\rm codim}}
\def\crd{\dim_{{\rm CR}}}
\def\crc{{\rm codim_{CR}}}
\def\phi{\varphi}
\def\eps{\varepsilon}
\def\d{\partial}
\def\a{\alpha}
\def\b{\beta}
\def\g{\gamma}
\def\G{\Gamma}
\def\D{\Delta}
\def\Om{\Omega}
\def\k{\kappa}
\def\l{\lambda}
\def\L{\Lambda}
\def\z{{\bar z}}
\def\w{{\bar w}}
\def\Z{{\mathbb Z}}
\def\t{\tau}
\def\th{\theta}

\emergencystretch15pt
\frenchspacing

\newtheorem{Thm}{Theorem}[section]
\newtheorem{Cor}[Thm]{Corollary}
\newtheorem{Pro}[Thm]{Proposition}
\newtheorem{Lem}[Thm]{Lemma}
\newtheorem{Prob}[Thm]{Problem}

\theoremstyle{definition}\newtheorem{Def}[Thm]{Definition}

\theoremstyle{remark}
\newtheorem{Rem}[Thm]{Remark}
\newtheorem{Exa}[Thm]{Example}
\newtheorem{Exs}[Thm]{Examples}

\def\bl{\begin{Lem}}
\def\el{\end{Lem}}
\def\bp{\begin{Pro}}
\def\ep{\end{Pro}}
\def\bt{\begin{Thm}}
\def\et{\end{Thm}}
\def\bc{\begin{Cor}}
\def\ec{\end{Cor}}
\def\bd{\begin{Def}}
\def\ed{\end{Def}}
\def\br{\begin{Rem}}
\def\er{\end{Rem}}
\def\be{\begin{Exa}}
\def\ee{\end{Exa}}
\def\bpr{\begin{Prob}}
\def\epr{\end{Prob}}
\def\bpf{\begin{proof}}
\def\epf{\end{proof}}
\def\ben{\begin{enumerate}} 
\def\een{\end{enumerate}}
\def\beq{\begin{equation}}
\def\eeq{\end{equation}}

\newcommand{ \ssf}{F^*}
\newcommand{ \jj}{\iota}
\newcommand{\alb}{\bar a_l}
\newcommand{\pz}{P(z,\bar z)}
\newcommand{\bet}{\beta}
\newcommand{\ga}{\gamma}
\newcommand{\de}{\delta}
\newcommand{\ka}{\kappa}
\newcommand{\ab}[1]{\vert z\vert^{#1}}
\newcommand{\cdva}{{\mathbb C^2}}
\newcommand{\om}{\Omega}
\newcommand{\ml}{M^{k,l}_a}
\newcommand{\te}{\theta}
\newcommand{\kl}{\dfrac k{l^2-k}}
\newcommand{\klv}{\dfrac k{\vert l^2-k\vert}}
\newcommand{\kldlo}{\dfrac{(4k-l^2-4)k^2}{(4k-4)(k^2-l^2)}}
\newcommand{ \ffn}{F^{\nu}}
\newcommand{ \fn}{f^{\nu}}
\newcommand{\pl}{P^{k,l}_a}
\newcommand{ \ksi}{\xi}
\newcommand{\zz}{(z,\bar z)}
\newcommand{ \nl}{\newline}
\newcommand{\fz}{{F(z,\bar z,u)}}
\newcommand{\fzs}{{\textstyle{F^*(z^*,\bar z^\vert (\alpha,\hat \al) \vert_{\La}*,u^*)}}}
\newcommand{\dom}{{\textstyle{b\Omega \ }}}
\newcommand{ \al}{\alpha}
\newcommand{\ppm}{m'}
\newcommand{\db}{\overline{\partial}}
\newcommand{\tf}{\tilde F}
\newcommand{\la}{\lambda}
\newcommand{\La}{\Lambda}
\newcommand{\kkk}{\frac{k}2 \frac{k}2}

\title[Characterization of  real-analytic infinitesimal CR automorphisms for a class of hypersurfaces in $\Bbb C^4.$
  ]{Characterization of   real-analytic  infinitesimal CR automorphisms for a class of hypersurfaces in $\Bbb C^4.$
}
\author{Cyril Julien and  Francine Meylan}
\address{ C. Julien: Department of Mathematics,
University of Fribourg, CH 1700 Perolles, Fribourg}

\email{cyril.julien@unifr.ch}
\address{ F. Meylan: Department of Mathematics,
University of Fribourg, CH 1700 Perolles, Fribourg}

\email{francine.meylan@unifr.ch}


\maketitle
\begin{abstract} 
In this paper,  motivated by the work of Kim and Kolar in  \cite{KK}  for the case of pseudoconvex models which are  sums of squares of polynomials, we study the Lie  algebra  $\5g$  of real-analytic  infinitesimal  $CR$ automorphisms of a  model hypersurface $M_0$ given by \begin{equation}\label{jm13}
 M_0= \{(z,w) \in \mathbb C^{3} \times  \mathbb C \ | \  \Im w= P\bar Q +  Q\bar P +  R\bar R  \},
\end{equation} where $P,$ $Q$ and $R$ are homogeneous polynomials.  In particular, we  classify $M_0$ with respect to  the  description of its  nilpotent rotations when $P,$ $Q$ and $R$ are monomials. We  also  give an example of a model $M_0$ for which  the real dimension of its  generalized (exotic) rotations is $3.$
\end{abstract}
\section{Introduction}
The structure of the  Lie algebra  $\5g$  of real-analytic  infinitesimal  $CR$ automorphisms of a  model hypersurface $M_0$  plays  a fundamental role when studying  the stability group of   a  hypersurface of finite type  in $\Bbb C^{n+1}$  whose model hypersurface is $M_0.$ This study  was initiated by Poincar\' e   for the unit sphere in $\Bbb C^2 $ \cite{Po} and pursued  by many others: see for instance   \cite{Be}, \cite{BFG}, \cite{C2}, \cite{IK}, \cite{KL},  \cite{Kr20}, \cite{Tanaka}, \cite{Tr}, \cite{V}, \cite{W}, \cite{Y}. 
 While Chern and Moser deal with the Levi non degenerate case in  \cite{CM}, that is, when the model hypersurface is an hyperquadric,  Kolar, Meylan and Zaitsev, inspired by the technics of Chern and Moser, deal with the Levi degenerate case in  \cite{KMZ}, that is, when the model hypersuface is  polynomial.
 
 \noindent More precisely,   Theorem 1.1 in \cite{KMZ} describes explicitely  the graded components  involved in the grading of $\5g,$  but does not indicate anything about the dimension or the vanishing of each graded component   regarding the model hypersurface.   

\noindent Recently, Kim and Kolar gave  in \cite{KK} a complete description of $\5g$ in the case of a pseudoconvex  model $M_0 \subset  \Bbb C^{n+1}$ which is a sum of squares, that is of the form $$\Im w = \sum_{j=1}^k P_j \overline{P_j},$$  where $P_j$ are homogeneous holomorphic polynomials,  under the hypothesis that $M_0$ is holomorphically nondegenerate (which implies  $k\ge n$).  In the case $k=n,$ it can be written using  the usual scalar product  as $$\Im w = <P,P>,$$ where  $ P=(P_1, \dots,P_n).$ This    motivates the following generalization:

\noindent   
 Let   $D$   be the symmetric matrix with respect to   the canonical basis $\{e_1, \dots ,e_n\}$ of $\Bbb C^n$  associated to the linear transformation 
sending $ e_j$ to $e_{\sigma(j)},$  where  $\sigma \in  S_n$  is a permutation that is a product of disjoint $2$-cycles. Let $P=(P_1, \dots, P_n)$ be a  $n$-tuple of  homogeneous  holomorphic  polynomials in the variables  $z=(z_1, \dots, z_n)$ for which  $^t\bar P D P $ is  a  (real)  homogeneous  polynomial.
We consider the model hypersurface denoted by  $M_d$ in $\Bbb C^{n+1}$   given by 
\begin{equation}\label{jm12}
\Im w= ^t\bar P D P , \ (z,w) \in \Bbb C^{n}\times \Bbb C.
\end{equation}
 Assuming that $M_d$ is holomorphically nondegenerate (which implies that  $P_j, \ j=1\dots, n,$ are linearly independent  over $\Bbb C$),  we address  the question of  characterizing  in terms of P  the Lie algebra $\5g$
of $M_d.$   

 \noindent 
 If the matrix $D$ is diagonal, that is, $M_d$ is pseudoconvex and  a sum of squares,   Kim and Kolar show that there are exactly  $3$ non vanishing  graded components, with only one depending on the polynomial $P_j.$ See Theorem 5.3 in \cite{KK}.

 \noindent 
 The case  $n=2$  has been solved   in \cite{FM17} and  \cite{KK}. 
 
 \noindent In this paper, we study  the case $n=3$  for non diagonal matrices $D,$ and we assume that  $P_1, P_2$ and $ P_3 $  are  monomial. We call this problem the {\it monomial  $PQR$ problem}. This  can also  be seen as the generalization  of the case $n=1$ given in  \cite{Ko1a}, \cite{KoMe}  where $P_1$ is necessarily a monomial of the form  $cz^k,$  and  $M_d$ in $\Bbb C^{2}$  is   given by 
\begin{equation}\label{jm012}
\Im w= |c|^2 |z|^{2k}, \ (z,w) \in \Bbb C^{}\times \Bbb C.
\end{equation} 
Note that $M_d$ given by \eqref{jm012} is the only Levi degenerate  model hypersurface  in  $\Bbb C^{2}$ that possesses a non linear symmetry.

 We will be particularly  interested into  the structure of the rigid  symmetries of weighted degree $0$, called the rotations. We will see that any such element can be decomposed into a diagonal rotation and a sum of at most $2$ nilpotent rotations. Note that the analogous result for $n=2$ is that any rotation is decomposed into the sum of a diagonal rotation and at most $1$ nilpotent rotation. Note that if $D$ is diagonal, the only rotations are the pure imaginary diagonal rotations (in Jordan form) \cite{KK}.
 
\noindent Another phenomenon we will discuss is the structure of the generalized rotations that impact the number of derivatives needed to determine completely an automorphism of a smooth perturbation of $M_0.$ Note that if $D$ is diagonal, there is no  such symmetries by Proposition 4.1 and Theorem 5.3 in \cite{KK}.

\noindent We emphasize that the study of  the structure of rotations may lead to the understanding (and hence the vanishing) of the graded components  of generalized rotations.  See   Proposition 4.1 in \cite{KK} and Section 3.

\noindent In Section 2, we recall  some definitions needed for the sequel. In Section 3, we collect, restate and  prove  some (partially) known results  on the graded  component of  $\5g$ of weight $0.$ In Section 4, 5 and 6, we study the graded components of  $\5g$ in the case of the monomial $PQR$ problem. Finally in Section 7, we state and prove the main results of the paper, Theorem\eqref{MainResult} and Theorem \eqref{MainResult2}.
\section{Preliminaries}

In this section, we recall the necessary definitions needed in the sequel regarding smooth hypersurfaces in $\mathbb C^{N}.$ 

Let  $M \subseteq \mathbb C^{n+1}$ be a  smooth hypersurface,
and $p \in M $ be a  point of {\em finite type} $d$ \cite{BG}.
We will consider
local holomorphic coordinates $(z,w)$ vanishing at $p$,
where $z =(z_1, z_2, ..., z_n)$ and  $z_j = x_j + iy_j$,
$w=u+iv$. The hyperplane $\{ v=0 \}$ is assumed to be tangent to
$M$ at $p$, hence  $M$  is described near $p$ as the graph of a uniquely
determined real valued function
\begin{equation} v = \psi(z_1,\dots, z_n,  \bar z_1,\dots,\bar z_n,  u), \ d\psi(0) = 0.
\label{vp1}
\end{equation}
{
  We may assume (\cite{BG}, \cite{K})
that
\begin{equation}\label{fifi}
\psi(z_1,\dots, z_n,  \bar z_1,\dots,\bar z_n,  u)=Q(z, \bar z) +o(|u|+|z|^d),
\end{equation}
where $Q(z, \bar z)$ is a nonzero homogeneous polynomial of degree $d$ with {no pluriharmonic} terms. Note that removing the pluriharmonic terms requires a change of coordinates }{in the variable $w$, which "absorbs`` the pluriharmonic terms into $w$.}

 { Reflecting} the  choice
of  variables given by  {\eqref{fifi}}, the
variables $w$, $u$ and $v$ are given weight $1,$  $\dfrac {\partial }{\partial {w}}$ weight $-1,$ while 
the complex tangential variables $(z_1, \dots, z_n)$  are given weight $ \dfrac{1}{d}$ and  $\dfrac {\partial }{\partial {z_{j}}}, \ j=1, \dots n, $ weight $ -\dfrac{1}{d}.$ 

\noindent  For this choice of weights, 
 {$Q(z, \bar z)$} is  a  weighted homogeneous polynomial of weighted  degree one, while $o(|u|+|z|^d)$ contains weighted homogeneous polynomials  of weighted degree  strictly bigger than one.

 Recall that the   derivatives of  weighted order $\kappa$  are  those of the form
 $$\dfrac{\partial^{\vert \alpha \vert + \vert \hat \alpha \vert + l}}{\partial z^{\alpha}\partial \bar z^{\hat \alpha} \partial u^l}, \ \ \ 
 \kappa=l +\dfrac{1}{d}\sum_{j=1}^n (\alpha_j + \hat \alpha_j).$$\emph{}

We may then rewrite  equation \eqref{fifi} as 
\begin{equation} v = 
\psi(z,  \bar z,  u)= 
Q\zz + o(1),
\label{vp12}
\end{equation}
where $o(1)$ denotes a smooth function  whose derivatives  at $0$ of weighted order less than or equal to
one vanish. This motivates the following definition 
\bd\cite {FM}{ Let $M$ be given by  (\ref{vp12}).
We define  {its associated
 model hypersurface 
by
\begin{equation} M_0: = \{(z,w) \in \mathbb C^{n+1}\ | \
 v  = Q \zz \}. \label{2}\end{equation}} }
We say that $M$ is a smooth perturbation of $M_Q.$
\ed
When  $Q= ^t\bar P D P,$  the associated model hypersurface is  given by  (\ref{jm12}), we write
\begin{equation}\label{jm13}
  M_d = \{(z,w) \in \mathbb C^{n+1}\ | \  \Im w= ^t\bar P D P   \}. 
\end{equation}
{\begin{Rem}
We may also consider  a $C^{k+d}$ perturbation  of $M_0, $ that is,  with  $C^{k+d}$ smooth function $o(1).$ 
\end{Rem}}
We denote  by $\Aut (M,p),$  the stability group of $M,$ that is, those germs at $p$
of biholomorphisms mapping $M$ into itself and fixing $p,$
and  by $\hol (M,p),$ the set of  germs of holomorphic vector fields in
$\mathbb C^{n+1}$ whose real part 
is tangent to $M$.
{
\begin{Rem}
If $M$ admits a holomorphic vector field $X$ in $\hol (M,p)$ such that $\Im X $ is also tangent,
then $\hol (M,p)$ is of infinite
 dimension (\cite{S}).
\end{Rem}}
We  write 
$$ \hol(M_0,p) = \bigoplus_{ \mu \ge -1 }  \5g_{\mu},$$
where $\5g_{\mu}$ consists of weighted homogeneous vector fields
of weight $\mu.$
 We recall the following definition.

\bd{A real-analytic hypersurface $M \subset \mathbb C^{n+1}$  is  {\it holomorphically nondegenerate
at $p \in M$} if there is no germ at $p$ of a holomorphic vector field $X$ tangent to $M.$}
 \ed

\section{The structure of $\5g_{0}$ and its consequences.}\label{ri}
We start  this section with a lemma that has been used several times (see for instance Proposition 4.1  in \cite{KK}) but whose proof was never given:  this lemma deals with the properties of a linear vector field  given in Jordan normal form. For the sake of clarity, we  give a proof of it.

\noindent Without loss of generality, we will  assume  that   the  nilpotent part  is  of the form  $\sum_{j \in J_1}  z_j  \dfrac {\partial }{\partial {z_{j+1}}}$,   where $J_1$ is the  subset of indices corresponding to the Jordan decomposition of the given vector field.
We recall the following definition 
\bd
 A linear vector field   in the variables $z$ given in Jordan normal form  is real diagonal  (respectively imaginary) if it  is diagonal with all the entries real (respectively imaginary). Without loss of generality, we say that $X$ is real (respectively imaginary)  diagonal if its  Jordan normal form is real (respectively imaginary) diagonal.
\ed
\bl  \label{great1} Let   $Q$  be a real homogeneous polynomial in $(z, \bar z)$  without pluriharmonic terms, and let 
 $W=X+Z$ be  a  holomorphic linear vector field  in the variables $z$ given  in Jordan normal form
with $X$ the diagonal and $Z$ the nilpotent part.
Then 
{
\begin{equation}
  W(Q)=Q \ \ \text {if and only if} \ \ \ X(Q)=Q,\  Z(Q)=0.
\end{equation}}
Moreover, if  the model hypersurface  given 
by
\begin{equation} M_Q: = \{(z,w) \in \mathbb C^{n+1}\ | \
 v  = Q \zz \}. \label{02}\end{equation} is holomorphically nondegenerate, then  $X$ is real diagonal and $Z=0.$
\el

\bpf  
Since $X$ is diagonal, we may consider  the spectral decomposition
$$Q=\sum Q_{\mu}, \quad Q_{\mu}\in \6Q_{\mu}:= \langle z^{\a}\bar z^{\b}:  X ( z^{\a}\bar z^{\b})= \mu z^{\a}\bar z^{\b}\rangle.$$
{Since   $[X,Z]=0,$    we obtain that $ Z(\6Q_{\mu})\subset \6Q_{\mu}$. }
We  put the lexicographic order on  monomials  $cz^{\a}\bar z^{\b}$ so that $ Z(z^{\a}\bar z^{\b})$  is   a sum of    monomials larger than $cz^{\a}\bar z^{\b}.$ 
We claim that $Q_{\mu}=0$ unless $\mu=1$.
Indeed, assume by contradiction that $Q_{\mu}\ne0$ for some $\mu\ne1$
and consider the minimal nontrivial monomial $c z^{\a}\bar z^{\b}$
in the expansion of $Q_{\mu}$ with respect to the lexicographic order.
Then, since $ Z(\6Q_{\mu})\subset \6Q_{\mu},$ $ Z(z^{\a}\bar z^{\b})$  is  a sum of   monomials   larger  than  $cz^{\a}\bar z^{\b}$ in $\6Q_{\mu}.$
Therefore the coefficient of $z^{\a}\bar z^{\b}$ in
$ W (Q)$ is equal to $\mu$. Since $ W (Q)=Q$, we must have $\mu=1$
contradicting our assumption.
Hence $Q=Q_{1}$ as claimed, implying $ X(Q)=Q$ and hence $Z(Q)=0.$
\epf
Let $X \in \hol(M_Q,p) $ be
given by 
$
X = \sum_{j=1}^n F_j(z,w) \partial_{ z_j} + G(z,w)\partial_{w}.
$
{We recall  that   $X$  is  rigid if $F_1, \dots, F_n, G $ are all independent of the variable  $w$  \cite{FM}.} The rigid  vector fields  of weighted degree 0  contained in $\5g_0$ are called rotations, while   $\5g_c$  denotes   the graded components  of  rigid  vector fields of weight  strictly bigger than $0,$ that is  with nonlinear coefficients; they are   called   generalized rotations or exotic symmetries. It is known that the weighted degree of an element of  $\5g_c$ is strictly less than one \cite{KMZ}.

\bt \label{great01} Let $M_0$ be given by  \eqref{2} and holomorphically nondegenerate.  Suppose   that $\5g_{0}$ contains no  imaginary  diagonal rotations. Then $\5g_1=0.$
\et
\begin{proof}
Suppose by contradiction that  $\5g_{1} \ne 0.$  Then, using a  result of \cite{KMZ}, we obtain that  $\dim \5g_{1}=1,$ with a generator  given by 
\begin{equation}\label{great002}
Y:=\sum \phi_j (z) w \partial _{z_j} + \dfrac{1}{2} w^2 \partial _{w}, \  \ 
2\sum \phi_j (z)  \partial _{z_j} (Q)=Q. 
\end{equation}
Let $S: =2\sum \phi_j (z)  \partial _{z_j}$ be the linear vector field associated to $Y.$ By Lemma \eqref{great1}, we may find linear holomorphic coordinates in  the $z$ variables  such that $S$ is real diagonal. We obtain then that   $iS$ is an imaginary diagonal  rotation since  $iS(Q)=iQ,$  and hence a contradiction.
\end{proof}
Since the grading element is in  $\5g_0, $  we have 
\bc  Let $M_0$ be given by  \eqref{2} and holomorphically nondegenerate. Suppose that  $\5g_{1} \ne 0.$  Then $\dim \5g_{0}>1.$
\ec
The following example shows that the converse is not true, unlike the case $n=2,$ \cite{FM17}.
\be
Consider the hypersurface  in  $\mathbb{C}^4$  given by 
	$$\Im w=
	(z_1^2z_2z_3^3
	+z_1z_2^4z_3)\overline{z_1}^2
	+{z_1}^2(\overline{z_1}^2\overline{z_2}\hspace{0.15mm}\overline{z_3}^3 +
	\overline{z_1}\hspace{0.15mm}\overline{z_2}^4\overline{z_3})
	+z_1^2z_2^2\overline{z_1}^2\overline{z_2}^2.$$
	The reader can check that $X$ given by 
	$$
	X:=i\hspace{.25mm}11z_1\frac{\partial}{\partial z_1}
	+i\hspace{.25mm}3z_2\frac{\partial}{\partial z_2}
	-i\hspace{.25mm}z_3\frac{\partial}{\partial z_3}.
	$$
	is a imaginary diagonal rotation. It is easy to check that $\dim \5g_{1}=0,$ using the characterization \eqref{great002}.

\ee
We now state     a theorem that has been proved in  \cite{KK} when  $M_0$  is pseudoconvex  (Proposition 4.1 in  \cite{KK}).  It shows the link between  the structure of the rotations and the generalized rotations, using a theorem of  \cite{KMZ}.  The proof is similar to the proof of Proposition 4.1 in  \cite{KK}.

\bt \label{great011} Let $M_0$ be given by  \eqref{2} be holomorphically nondegenerate. Suppose that $\5g_{1} \ne 0$ and that $\5g_{0}$ contains no  real diagonal rotations. Then $\5g_c=0.$
\et
\begin{proof}
Since  $\5g_{1} \ne 0,$ using a  result of \cite{KMZ}, we obtain that  $\dim \5g_{1}=1,$ with a generator  given by 
\begin{equation}\label{great02}
Y:= \sum \phi_j (z) w \partial _{z_j} + \dfrac{1}{2} w^2 \partial _{w}, \  \ 
2\sum \phi_j (z)  \partial _{z_j} (Q)=Q. 
\end{equation}
Let $S: =2\sum \phi_j (z)  \partial _{z_j}$ be the linear vector field associated to $Y.$ By Lemma \eqref{great1}, we may find linear holomorphic coordinates in  the $z$ variables  such that $S$ is real diagonal. Without loss of generality, we  may then assume that $Y$ is  of the form
$$Y= \sum \lambda_j z_j w \partial _{z_j} + \dfrac{1}{2} w^2 \partial _{w},\ \ \lambda_j \in  \Bbb R.$$
Assume by contradiction that there exists  a non zero $X \in \5g_c$  of weight  $\mu>0.$ Then again by a result of \cite{KMZ}, the Lie bracket $[X, Y]=0,$  the reason beeing that there is no graded component of weight strictly bigger than one under the assumption  that $M_0$ is holomorphiacally nondegenerate.
Let $E  \in \5g_{0}$ be the 
grading element given by 
\begin{equation}\label{gentil.02}
 E = \sum_{j=1}^n \frac{1}{d} z_j \partial_{ z_j} + w\partial_{w}.
\end{equation}
Since  $[E, X] =\mu X ,$   $[\partial_ w, X]=0$ and  $[ [\partial_ w, Y], X]=0$ by the identity of Jacobi,  we obtain that
$$E-[\partial_ w, Y]= \sum (\frac{1}{d}-\lambda_j )z_j  \partial _{z_j}$$ is a nonzero real diagonal  rotation, and hence a contradiction.

\end{proof}

\section{ The structure of $\5g_{0}$ for  the monomial $PQR$ problem }
In this section, we consider the  monomial $PQR$ problem, that is, when $M_d$ is   given by 
\begin{equation}\label{jm13}
 M_0= \{(z,w) \in \mathbb C^{3} \times  \mathbb C \ | \  \Im w= P\bar Q +  Q\bar P +  R\bar R  \},
\end{equation}
where $P,$ $Q$ and $R$ are monomials. We show  that any rotation of $M_d$  can be decomposed  into a  rotation  that is diagonal {\it in these same coordinates}  and a sum of  at most  $2$ nilpotent rotations. We start with the following  lemma  that gives a characterization of holomorphic nondegeneracy.
\bl\label{LOG1}
	Let $M_d$ be given by \eqref {jm13}, with $P,$ $Q$ and $R$  not necessarily monomials.
Then  $M_d$ is holomorphically nondegenerate if and only if  the  jacobian  of  $P$, $Q$ and  $R$ is nonzero.
\el

\begin{proof} Let $J(P,Q,R)$ be the jacobian of $P,$ $Q$ and $R.$   By contradiction, if      $J(P,Q,R) =0,$  solving a system of $3$ equations, one can find  a holomorphic vector field  
 	$
	X = \sum_{j=1}^{3}f_j(z)\frac{\partial}{\partial z_j} 
	$
	satisfying 
	$
	X(P)=0,$ $X(Q)=0$ and $X(R)=0.$
	  Hence the contradiction.

	Suppose now that  $J(P,Q,R)\neq 0.$ 
	Let $X$ be a holomorphic vector field  tangent to $M_d$ given by 
	\begin{equation}
	X=\sum_{j=1}^{3}f_j(z,w)\frac{\partial}{\partial z_j} +  g(z,w)\frac{\partial}{\partial w}.
	\end{equation}
	This implies   $g(z,u) =0,$ and hence 
	\begin{equation} g(z,w) =0, \ 
	X(P)\overline{Q}+X(Q)\overline{P}+X(R)\overline{R}=0,
	\end{equation}
	and hence, by  assumption, 
	$$X(P)=0,\
	X(Q)=0,\
	X(R)=0.$$ This yields to a Cramer  system of $3$ equations whose only solution is the zero solution since  $J(P,Q,R)\neq 0.$  Hence $X=0,$ that is,  $M_d$ is holomorphically nondegenerate
	
	\end{proof}

\bp~\label{LOG}

	Let  $M_d$ be given by \eqref {jm13}
	with 
	$P=z_1^{\alpha_1}z_2^{\alpha_2}z_3^{\alpha_3}$,
	$Q=z_1^{\beta_1}z_2^{\beta_2}z_3^{\beta_3}$ et
	$R=z_1^{\gamma_1}z_2^{\gamma_2}z_3^{\gamma_3}$.
	Then  $M$ is holomorphically nondegenerate  if and only if 
	$$
	\begin{vmatrix}
	\alpha_1&\beta_1&\gamma_1\\
	\alpha_2&\beta_2&\gamma_2\\
	\alpha_3&\beta_3&\gamma_3
	\end{vmatrix}
	\neq0.
	$$
\ep

\begin{proof} 
	By Lemma\eqref{LOG1}, we have $M_d$ holomorphically nondegenerate if and only if  
	
	\begin{eqnarray}
	0&\neq& J(P,Q,R)
	=
		\begin{vmatrix}
	P_{z_1}&Q_{z_1}&R_{z_1}\\
	P_{z_2}	&Q_{z_2}&	R_{z_2}\\
	P_{z_3}&	Q_{z_3}&	R_{z_3}
	\end{vmatrix}
	\nonumber\\\nonumber\\
	&=&
		\begin{vmatrix}
	\alpha_1z_1^{\alpha_1-1}z_2^{\alpha_2}z_3^{\alpha_3}
	&\beta_1z_1^{\beta_1-1}z_2^{\beta_2}z_3^{\beta_3}
	&\gamma_1z_1^{\gamma_1-1}z_2^{\gamma_2}z_3^{\gamma_3}
	\\
	\alpha_2z_1^{\alpha_1}z_2^{\alpha_2-1}z_3^{\alpha_3}
	&\beta_2z_1^{\beta_1}z_2^{\beta_2-1}z_3^{\beta_3}
	&\gamma_2z_1^{\gamma_1}z_2^{\gamma_2-1}z_3^{\gamma_3}
	\\
	\alpha_3z_1^{\alpha_1}z_2^{\alpha_2}z_3^{\alpha_3-1}
	&\beta_3z_1^{\beta_1}z_2^{\beta_2}z_3^{\beta_3-1}
	&\gamma_3z_1^{\gamma_1}z_2^{\gamma_2}z_3^{\gamma_3-1}
	\end{vmatrix}
	\nonumber
	\\\nonumber\\
	&=&
	\textcolor{white}{+}
	\alpha_1z_1^{\alpha_1-1}z_2^{\alpha_2}z_3^{\alpha_3}
		\begin{vmatrix}
	\beta_2z_1^{\beta_1}z_2^{\beta_2-1}z_3^{\beta_3}
	&\gamma_2z_1^{\gamma_1}z_2^{\gamma_2-1}z_3^{\gamma_3}
	\\
	\beta_3z_1^{\beta_1}z_2^{\beta_2}z_3^{\beta_3-1}
	&\gamma_3z_1^{\gamma_1}z_2^{\gamma_2}z_3^{\gamma_3-1}
	\end{vmatrix}
	\nonumber\\\nonumber\\
	&&-
	\beta_1z_1^{\beta_1-1}z_2^{\beta_2}z_3^{\beta_3}
			\begin{vmatrix}
	\alpha_2z_1^{\alpha_1}z_2^{\alpha_2-1}z_3^{\alpha_3}
	&\gamma_2z_1^{\gamma_1}z_2^{\gamma_2-1}z_3^{\gamma_3}
	\\
	\alpha_3z_1^{\alpha_1}z_2^{\alpha_2}z_3^{\alpha_3-1}
	&\gamma_3z_1^{\gamma_1}z_2^{\gamma_2}z_3^{\gamma_3-1}
	\end{vmatrix}
	\nonumber\\\nonumber\\
	&&
	+
	\gamma_1z_1^{\gamma_1-1}z_2^{\gamma_2}z_3^{\gamma_3}
			\begin{vmatrix}
	\alpha_2z_1^{\alpha_1}z_2^{\alpha_2-1}z_3^{\alpha_3}
	&\beta_2z_1^{\beta_1}z_2^{\beta_2-1}z_3^{\beta_3}
	\\
	\alpha_3z_1^{\alpha_1}z_2^{\alpha_2}z_3^{\alpha_3-1}
	&\beta_3z_1^{\beta_1}z_2^{\beta_2}z_3^{\beta_3-1}
	\end{vmatrix}
	\nonumber
	\\\nonumber\\
	&=&
		\begin{vmatrix}
	\alpha_1&\beta_1&\gamma_1\\
	\alpha_2&\beta_2&\gamma_2\\
	\alpha_3&\beta_3&\gamma_3
	\end{vmatrix}
z_1^{\alpha_1+\beta_1+\gamma_1-1}
z_2^{\alpha_2+\beta_2+\gamma_2-1}
z_3^{\alpha_3+\beta_3+\gamma_3-1}.
	\nonumber
	\end{eqnarray}

\end{proof}

We may now state the main result of this section.

\bt\label{L33}
	Let  $M_d$ be given by \eqref {jm13}
	with $P,$ $Q$ and $R$ monomials, and holomorphically nondegenerate.
	\noindent Then any $X \in \5g_{0} $ given by 
$$	X=\sum_{j,k=1}^{3}a_{jk}z_k\frac{\partial}{\partial z_j},\;\;\,a_{jk}\in\mathbb{C},
	$$ may be written as 
		$$
	X=D+N,
	$$ where \begin{equation}\label{lilia0}D = \sum_{j=1 }^{3}a_{jj}z_j\frac{\partial}{\partial z_j} \in \5g_{0},\end{equation} 
	 and \begin{equation}\label{lilia}N  = \sum_{j\ne k, j,k=1}^{3}a_{jk}z_k\frac{\partial}{\partial z_j} \in \5g_{0},\end{equation} satisfies the following property:
	 \begin{enumerate}
	 \item either $N$ is  a nilpotent rotation,  
	 \item or there exists  a subset of indices $J_1$ of $J =\{ (j,k), j\ne k, j,k=1, 2, 3,\}$ such that $$N_1: ={\sum}_{J_1} a_{jk}z_k\frac{\partial}{\partial z_j} \ \text{and} \  \  N_2: ={\sum}_{J-J_1}a_{jk}z_k\frac{\partial}{\partial z_j}$$  are nilpotent rotations, with $N=N_1 +N_2.$
	\end {enumerate}
\et
\begin{proof}
 We write 
$$P=c_Pz_1^{\alpha_1} z_2^{\alpha_2}z_3^{\alpha_3}=c_Pz^{\alpha} ,\ 
Q=c_Qz_1^{\beta_1} z_2^{\beta_2}z_3^{\beta_3}=c_Qz^{\beta},  $$
$$R=c_Rz_1^{\gamma_1} z_2^{\gamma_2}z_3^{\gamma_3} = c_Rz^{\gamma}, \ c_P,\  c_Q, \  c_R \in  \Bbb C.$$  
We first show that $D\in \5g_{0}.$ 

\noindent By assumption, we have $\Re X( P\bar Q +  Q\bar P +  R\bar R ) = \Re D( P\bar Q +  Q\bar P +  R\bar R ) + \Re N( P\bar Q +  Q\bar P +  R\bar R ) =0.$

\noindent Observing  that the terms given by $\Re N( P\bar Q +  Q\bar P +  R\bar R )$
are 
\scriptsize{
	\begin{equation}\label{stresa}
	\begin{array}{lll}
	&{\color{white}{+}}
	a_{12}c_P\overline{c_Q}\alpha_1
	z_1^{\alpha_1-1} z_2^{\alpha_2+1}z_3^{\alpha_3}
	\overline{z_1}^{\beta_1} \overline{z_2}^{\beta_2}\overline{z_3}^{\beta_3}
	&+
	a_{13}c_P\overline{c_Q}\alpha_1
	z_1^{\alpha_1-1} z_2^{\alpha_2}z_3^{\alpha_3+1}
	\overline{z_1}^{\beta_1} \overline{z_2}^{\beta_2}\overline{z_3}^{\beta_3}
	\\
	&+
	\overline{a_{12}}c_Q\overline{c_P}\alpha_1
	{z_1}^{\beta_1} {z_2}^{\beta_2}{z_3}^{\beta_3}
	\overline{z_1}^{\alpha_1-1} \overline{z_2}^{\alpha_2+1}\overline{z_3}^{\alpha_3}
	&+
	\overline{a_{13}}c_Q\overline{c_P}\alpha_1
	{z_1}^{\beta_1} {z_2}^{\beta_2}{z_3}^{\beta_3}
	\overline{z_1}^{\alpha_1-1} \overline{z_2}^{\alpha_2}\overline{z_3}^{\alpha_3+1}
	\\
	\\
	&+
	a_{12}c_Q\overline{c_P}\beta_1
	z_1^{\beta_1-1} z_2^{\beta_2+1}z_3^{\beta_3}
	\overline{z_1}^{\alpha_1} \overline{z_2}^{\alpha_2}\overline{z_3}^{\alpha_3}
	&+
	a_{13}c_Q\overline{c_P}\beta_1
	z_1^{\beta_1-1} z_2^{\beta_2}z_3^{\beta_3+1}
	\overline{z_1}^{\alpha_1} \overline{z_2}^{\alpha_2}\overline{z_3}^{\alpha_3}
	\\
	&+
	\overline{a_{12}}c_P\overline{c_Q}\beta_1
	{z_1}^{\alpha_1} {z_2}^{\alpha_2}{z_3}^{\alpha_3}
	\overline{z_1}^{\beta_1-1} \overline{z_2}^{\beta_2+1}\overline{z_3}^{\beta_3}
	&+
	\overline{a_{13}}c_P\overline{c_Q}\beta_1
	{z_1}^{\alpha_1} {z_2}^{\alpha_2}{z_3}^{\alpha_3}
	\overline{z_1}^{\beta_1-1} \overline{z_2}^{\beta_2}\overline{z_3}^{\beta_3+1}
	\\
	\\
	&+
	a_{12}c_R\overline{c_R}\gamma_1
	z_1^{\gamma_1-1} z_2^{\gamma_2+1}z_3^{\gamma_3}
	\overline{z_1}^{\gamma_1} \overline{z_2}^{\gamma_2}\overline{z_3}^{\gamma_3}
	&+
	a_{13}c_R\overline{c_R}\gamma_1
	z_1^{\gamma_1-1} z_2^{\gamma_2}z_3^{\gamma_3+1}
	\overline{z_1}^{\gamma_1} \overline{z_2}^{\gamma_2}\overline{z_3}^{\gamma_3}
	\\
	
	&+
	\overline{a_{12}}c_R\overline{c_R}\gamma_1
	{z_1}^{\gamma_1} {z_2}^{\gamma_2}{z_3}^{\gamma_3}
	\overline{z_1}^{\gamma_1-1} \overline{z_2}^{\gamma_2+1}\overline{z_3}^{\gamma_3}
	&+
	\overline{a_{13}}c_R\overline{c_R}\gamma_1
	{z_1}^{\gamma_1} {z_2}^{\gamma_2}{z_3}^{\gamma_3}
	\overline{z_1}^{\gamma_1-1} \overline{z_2}^{\gamma_2}\overline{z_3}^{\gamma_3+1}
	\\
	\\
	\\
	\\
	+a_{21}c_P\overline{c_Q}\alpha_2
	z_1^{\alpha_1+1} z_2^{\alpha_2-1}z_3^{\alpha_3}
	\overline{z_1}^{\beta_1} \overline{z_2}^{\beta_2}\overline{z_3}^{\beta_3}
	&
	&+
	a_{23}c_P\overline{c_Q}\alpha_2
	z_1^{\alpha_1} z_2^{\alpha_2-1}z_3^{\alpha_3+1}
	\overline{z_1}^{\beta_1} \overline{z_2}^{\beta_2}\overline{z_3}^{\beta_3}
	\\
	+
	\overline{a_{21}}c_Q\overline{c_P}\alpha_2
	{z_1}^{\beta_1} {z_2}^{\beta_2}{z_3}^{\beta_3}
	\overline{z_1}^{\alpha_1+1} \overline{z_2}^{\alpha_2-1}\overline{z_3}^{\alpha_3}
	&
	&+
	\overline{a_{23}}c_Q\overline{c_P}\alpha_2
	{z_1}^{\beta_1} {z_2}^{\beta_2}{z_3}^{\beta_3}
	\overline{z_1}^{\alpha_1} \overline{z_2}^{\alpha_2-1}\overline{z_3}^{\alpha_3+1}
	\\
	\\
	+a_{21}c_Q\overline{c_P}\beta_2
	z_1^{\beta_1+1} z_2^{\beta_2-1}z_3^{\beta_3}
	\overline{z_1}^{\alpha_1} \overline{z_2}^{\alpha_2}\overline{z_3}^{\alpha_3}
	&
	&+
	a_{23}c_Q\overline{c_P}\beta_2
	z_1^{\beta_1} z_2^{\beta_2-1}z_3^{\beta_3+1}
	\overline{z_1}^{\alpha_1} \overline{z_2}^{\alpha_2}\overline{z_3}^{\alpha_3}
	\\
	+
	\overline{a_{21}}c_P\overline{c_Q}\beta_2
	{z_1}^{\alpha_1} {z_2}^{\alpha_2}{z_3}^{\alpha_3}
	\overline{z_1}^{\beta_1+1} \overline{z_2}^{\beta_2-1}\overline{z_3}^{\beta_3}
	&
	&+
	\overline{a_{23}}c_P\overline{c_Q}\beta_2
	{z_1}^{\alpha_1} {z_2}^{\alpha_2}{z_3}^{\alpha_3}
	\overline{z_1}^{\beta_1} \overline{z_2}^{\beta_2-1}\overline{z_3}^{\beta_3+1}
	\\
	\\
	+a_{21}c_R\overline{c_R}\gamma_2
	z_1^{\gamma_1+1} z_2^{\gamma_2-1}z_3^{\gamma_3}
	\overline{z_1}^{\gamma_1} \overline{z_2}^{\gamma_2}\overline{z_3}^{\gamma_3}
	&
	&+
	a_{23}c_R\overline{c_R}\gamma_2
	z_1^{\gamma_1} z_2^{\gamma_2-1}z_3^{\gamma_3+1}
	\overline{z_1}^{\gamma_1} \overline{z_2}^{\gamma_2}\overline{z_3}^{\gamma_3}
	\\
	+
	\overline{a_{21}}c_R\overline{c_R}\gamma_2
	{z_1}^{\gamma_1} {z_2}^{\gamma_2}{z_3}^{\gamma_3}
	\overline{z_1}^{\gamma_1+1} \overline{z_2}^{\gamma_2-1}\overline{z_3}^{\gamma_3}
	&
	&+
	\overline{a_{23}}c_R\overline{c_R}\gamma_2
	{z_1}^{\gamma_1} {z_2}^{\gamma_2}{z_3}^{\gamma_3}
	\overline{z_1}^{\gamma_1} \overline{z_2}^{\gamma_2-1}\overline{z_3}^{\gamma_3+1}
	\\
	\\
	\\
	\\
	+a_{31}c_P\overline{c_Q}\alpha_3
	z_1^{\alpha_1+1} z_2^{\alpha_2}z_3^{\alpha_3-1}
	\overline{z_1}^{\beta_1} \overline{z_2}^{\beta_2}\overline{z_3}^{\beta_3}
	&+
	a_{32}c_P\overline{c_Q}\alpha_3
	z_1^{\alpha_1} z_2^{\alpha_2+1}z_3^{\alpha_3-1}
	\overline{z_1}^{\beta_1} \overline{z_2}^{\beta_2}\overline{z_3}^{\beta_3}
	&
	\\
	+
	\overline{a_{31}}c_Q\overline{c_P}\alpha_3
	{z_1}^{\beta_1} {z_2}^{\beta_2}{z_3}^{\beta_3}
	\overline{z_1}^{\alpha_1+1} \overline{z_2}^{\alpha_2}\overline{z_3}^{\alpha_3-1}
	&+
	\overline{a_{32}}c_Q\overline{c_P}\alpha_3
	{z_1}^{\beta_1} {z_2}^{\beta_2}{z_3}^{\beta_3}
	\overline{z_1}^{\alpha_1} \overline{z_2}^{\alpha_2+1}\overline{z_3}^{\alpha_3-1}
	&
	\\
	\\
	+a_{31}c_Q\overline{c_P}\beta_3
	z_1^{\beta_1+1} z_2^{\beta_2}z_3^{\beta_3-1}
	\overline{z_1}^{\alpha_1} \overline{z_2}^{\alpha_2}\overline{z_3}^{\alpha_3}
	&+
	a_{32}c_Q\overline{c_P}\beta_3
	z_1^{\beta_1} z_2^{\beta_2+1}z_3^{\beta_3-1}
	\overline{z_1}^{\alpha_1} \overline{z_2}^{\alpha_2}\overline{z_3}^{\alpha_3}
	&
	\\
	+
	\overline{a_{31}}c_P\overline{c_Q}\beta_3
	{z_1}^{\alpha_1} {z_2}^{\alpha_2}{z_3}^{\alpha_3}
	\overline{z_1}^{\beta_1+1} \overline{z_2}^{\beta_2}\overline{z_3}^{\beta_3-1}
	&+
	\overline{a_{32}}c_P\overline{c_Q}\beta_3
	{z_1}^{\alpha_1} {z_2}^{\alpha_2}{z_3}^{\alpha_3}
	\overline{z_1}^{\beta_1} \overline{z_2}^{\beta_2+1}\overline{z_3}^{\beta_3-1}
	&
	\\
	\\
	+a_{31}c_R\overline{c_R}\gamma_3
	z_1^{\gamma_1+1} z_2^{\gamma_2}z_3^{\gamma_3-1}
	\overline{z_1}^{\gamma_1} \overline{z_2}^{\gamma_2}\overline{z_3}^{\gamma_3}
	&+
	a_{32}c_R\overline{c_R}\gamma_3
	z_1^{\gamma_1} z_2^{\gamma_2+1}z_3^{\gamma_3-1}
	\overline{z_1}^{\gamma_1} \overline{z_2}^{\gamma_2}\overline{z_3}^{\gamma_3}
	&
	\\
	+
	\overline{a_{31}}c_R\overline{c_R}\gamma_3
	{z_1}^{\gamma_1} {z_2}^{\gamma_2}{z_3}^{\gamma_3}
	\overline{z_1}^{\gamma_1+1} \overline{z_2}^{\gamma_2}\overline{z_3}^{\gamma_3-1}
	&+
	\overline{a_{32}}c_R\overline{c_R}\gamma_3
	{z_1}^{\gamma_1} {z_2}^{\gamma_2}{z_3}^{\gamma_3}
	\overline{z_1}^{\gamma_1} \overline{z_2}^{\gamma_2+1}\overline{z_3}^{\gamma_3-1},
	&
	\end{array}
	\end{equation}
	}\normalsize
and using Proposition  \eqref{LOG}, we obtain that the terms of $\Re D(  P\bar Q +  Q\bar P +  R\bar R ),$  which are of the form
$z^{\alpha} 
\overline{z}^{\beta}, \  
z^{\beta} 
\overline{z}^{\alpha},  \  
z^{\gamma} 
\overline{z}^{\gamma},$
can never cancel out with the terms of \eqref{stresa}.

\noindent  Therefore  $\Re D( P\bar Q +  Q\bar P +  R\bar R ) =0,$ that is $D\in \5g_{0},$ and hence $N \in \5g_{0}.$ We need the following lemma.

\bl
\label{TheLemma}~
	
	\nopagebreak[4]\noindent
	Let 
	$$
	N=(a_{12}z_2+a_{13}z_3)\frac{\partial}{\partial z_1}
	+(a_{21}z_1+a_{23}z_3)\frac{\partial}{\partial z_2}
	+(a_{31}z_1+a_{32}z_2)\frac{\partial}{\partial z_3} \in \5g_{0}, \ a_{12}\neq0. \ 
$$    Then 
	\begin{enumerate}
	\item  either  $$N_1=a_{12}z_2\frac{\partial}{\partial z_1} \in \5g_{0},$$   and  $M_d $ is given by  
$	\Im w=
	2\Re\left(c_P\overline{c_Q}
	z_1z_2^{\alpha_2}z_3^{\alpha_3}\overline{z_2}^{\alpha_2+1}\overline{z_3}^{\alpha_3}
	\right)	
	+
	c_R\overline{c_R}
	z_2^{\gamma_2}z_3^{\gamma_3}\overline{z_2}^{\gamma_2}\overline{z_3}^{\gamma_3},$
	\\
	\item or  $$ N_1= a_{12}z_2\frac{\partial}{\partial z_1} + a_{23}z_3\frac{\partial}{\partial z_2}\in \5g_{0},$$  and  $M_d $ is given by 
	$\Im w=
	2\Re\left(c_P\overline{c_Q}
	z_1z_3^k\overline{z_3}^{k+1}
	\right)
	+
	c_R\overline{c_R}
	z_2z_3^k\overline{z_2}\hspace{.15mm}\overline{z_3}^k,
	$ \\
	\item or $$N_1= a_{12}z_2\frac{\partial}{\partial z_1} + a_{31}z_1\frac{\partial}{\partial z_3}\in \5g_{0},$$ and  $M_d $ is given by 
	$\Im w=
	2\Re\left(c_P\overline{c_Q}
	z_2^kz_3\overline{z_2}^{k+1}
	\right)
 +
	c_R\overline{c_R}
	z_1z_2^k\overline{z_1}\hspace{.15mm}\overline{z_2}^k. $ \\
	
	\end{enumerate}

\el
\begin{proof} Suppose that  $N_1=a_{12}z_2\frac{\partial}{\partial z_1} \in \5g_{0}.$
It implies that 

\begin{equation}
\begin{array}{ll}
0=
&{\color{white}{+}}
a_{12}c_P\overline{c_Q}\alpha_1
z_1^{\alpha_1-1} z_2^{\alpha_2+1}z_3^{\alpha_3}
\overline{z_1}^{\beta_1} \overline{z_2}^{\beta_2}\overline{z_3}^{\beta_3}
\\
&+
\overline{a_{12}}c_Q\overline{c_P}\alpha_1
{z_1}^{\beta_1} {z_2}^{\beta_2}{z_3}^{\beta_3}
\overline{z_1}^{\alpha_1-1} \overline{z_2}^{\alpha_2+1}\overline{z_3}^{\alpha_3}
\\
\\
&+
a_{12}c_Q\overline{c_P}\beta_1
z_1^{\beta_1-1} z_2^{\beta_2+1}z_3^{\beta_3}
\overline{z_1}^{\alpha_1} \overline{z_2}^{\alpha_2}\overline{z_3}^{\alpha_3}
\\
&+
\overline{a_{12}}c_P\overline{c_Q}\beta_1
{z_1}^{\alpha_1} {z_2}^{\alpha_2}{z_3}^{\alpha_3}
\overline{z_1}^{\beta_1-1} \overline{z_2}^{\beta_2+1}\overline{z_3}^{\beta_3}
\\
\\
&+
a_{12}c_R\overline{c_R}\gamma_1
z_1^{\gamma_1-1} z_2^{\gamma_2+1}z_3^{\gamma_3}
\overline{z_1}^{\gamma_1} \overline{z_2}^{\gamma_2}\overline{z_3}^{\gamma_3}
\\
&+
\overline{a_{12}}c_R\overline{c_R}\gamma_1
{z_1}^{\gamma_1} {z_2}^{\gamma_2}{z_3}^{\gamma_3}
\overline{z_1}^{\gamma_1-1} \overline{z_2}^{\gamma_2+1}\overline{z_3}^{\gamma_3} .
\end{array}
\end{equation}
 By Proposition\eqref{LOG}, there  are $2$ possibilities: 
 \begin{enumerate}
 \item
 either 
$$
a_{12}c_P\overline{c_Q}\alpha_1
z_1^{\alpha_1-1} z_2^{\alpha_2+1}z_3^{\alpha_3}
\overline{z_1}^{\beta_1} \overline{z_2}^{\beta_2}\overline{z_3}^{\beta_3}
\;\;\;\;
\text{matches  with}
\;\;\;\;
\overline{a_{12}}c_Q\overline{c_P}\alpha_1
{z_1}^{\beta_1} {z_2}^{\beta_2}{z_3}^{\beta_3}
\overline{z_1}^{\alpha_1-1} \overline{z_2}^{\alpha_2+1}\overline{z_3}^{\alpha_3}
$$ \\
\item 
or
$$
a_{12}c_P\overline{c_Q}\alpha_1
z_1^{\alpha_1-1} z_2^{\alpha_2+1}z_3^{\alpha_3}
\overline{z_1}^{\beta_1} \overline{z_2}^{\beta_2}\overline{z_3}^{\beta_3}
\;\;\;\;
\text{matches with}
\;\;\;\;
\overline{a_{12}}c_R\overline{c_R}\gamma_1
{z_1}^{\gamma_1} {z_2}^{\gamma_2}{z_3}^{\gamma_3}
\overline{z_1}^{\gamma_1-1} \overline{z_2}^{\gamma_2+1}\overline{z_3}^{\gamma_3}.
$$
\end{enumerate}
Case (1) yields to 
$$
\begin{array}{ccc}
\alpha_1-1&=&\beta_1\\
\alpha_2+1&=&\beta_2\\
\alpha_3&=&\beta_3,
\end{array}
$$
which implies 
$$
\begin{array}{ccccc}
\gamma_1=0&&\alpha_1=1&&\beta_1=0.\\
\end{array}
$$
Hence $M_d$ is given by 
\begin{equation}
\Im w=
c_P\overline{c_Q}
z_1z_2^{\alpha_2}z_3^{\alpha_3}\overline{z_2}^{\alpha_2+1}\overline{z_3}^{\alpha_3}
+
c_Q\overline{c_P}
{z_2}^{\alpha_2+1}{z_3}^{\alpha_3}
\overline{z_1}\hspace{.15mm}\overline{z_2}^{\alpha_2}\overline{z_3}^{\alpha_3}
+
c_R\overline{c_R}
{z_2}^{\gamma_2}{z_3}^{\gamma_3}
\overline{z_2}^{\gamma_2}\overline{z_3}^{\gamma_3}.
\end{equation}

Case (2) yields to 

\begin{equation}\label{idr}
\begin{array}{ccccc}
\gamma_1&=&\alpha_1-1&=&\beta_1+1\\
\gamma_2&=&\alpha_2+1&=&\beta_2-1\\
\gamma_3&=&\alpha_3&=&\beta_3.
\end{array}
\end{equation}
Using Proposition\eqref{LOG}, \eqref{idr} implies that $M_d$ is holomorphically degenerate: case (2) does not occur.

\noindent  For the rest of the proof, we assume  that $N_1=a_{12}z_2\frac{\partial}{\partial z_1} \not \in \5g_{0}$ and suppose  $a_{13}\neq0.$

\noindent We  claim  that  terms containing  $a_{12}$ can  not interfer with terms containing  $a_{13}.$ 

\noindent The terms under study are 
\begin{equation} 
\begin{array}{lll}
&{\color{white}{+}}
a_{12}c_P\overline{c_Q}\alpha_1
z_1^{\alpha_1-1} z_2^{\alpha_2+1}z_3^{\alpha_3}
\overline{z_1}^{\beta_1} \overline{z_2}^{\beta_2}\overline{z_3}^{\beta_3}
&+
a_{13}c_P\overline{c_Q}\alpha_1
z_1^{\alpha_1-1} z_2^{\alpha_2}z_3^{\alpha_3+1}
\overline{z_1}^{\beta_1} \overline{z_2}^{\beta_2}\overline{z_3}^{\beta_3}
\\
&+
\overline{a_{12}}c_Q\overline{c_P}\alpha_1
{z_1}^{\beta_1} {z_2}^{\beta_2}{z_3}^{\beta_3}
\overline{z_1}^{\alpha_1-1} \overline{z_2}^{\alpha_2+1}\overline{z_3}^{\alpha_3}
&+
\overline{a_{13}}c_Q\overline{c_P}\alpha_1
{z_1}^{\beta_1} {z_2}^{\beta_2}{z_3}^{\beta_3}
\overline{z_1}^{\alpha_1-1} \overline{z_2}^{\alpha_2}\overline{z_3}^{\alpha_3+1}
\\
\\
&+
a_{12}c_Q\overline{c_P}\beta_1
z_1^{\beta_1-1} z_2^{\beta_2+1}z_3^{\beta_3}
\overline{z_1}^{\alpha_1} \overline{z_2}^{\alpha_2}\overline{z_3}^{\alpha_3}
&+
a_{13}c_Q\overline{c_P}\beta_1
z_1^{\beta_1-1} z_2^{\beta_2}z_3^{\beta_3+1}
\overline{z_1}^{\alpha_1} \overline{z_2}^{\alpha_2}\overline{z_3}^{\alpha_3}
\\
&+
\overline{a_{12}}c_P\overline{c_Q}\beta_1
{z_1}^{\alpha_1} {z_2}^{\alpha_2}{z_3}^{\alpha_3}
\overline{z_1}^{\beta_1-1} \overline{z_2}^{\beta_2+1}\overline{z_3}^{\beta_3}
&+
\overline{a_{13}}c_P\overline{c_Q}\beta_1
{z_1}^{\alpha_1} {z_2}^{\alpha_2}{z_3}^{\alpha_3}
\overline{z_1}^{\beta_1-1} \overline{z_2}^{\beta_2}\overline{z_3}^{\beta_3+1}
\\
\\
&+
a_{12}c_R\overline{c_R}\gamma_1
z_1^{\gamma_1-1} z_2^{\gamma_2+1}z_3^{\gamma_3}
\overline{z_1}^{\gamma_1} \overline{z_2}^{\gamma_2}\overline{z_3}^{\gamma_3}
&+
a_{13}c_R\overline{c_R}\gamma_1
z_1^{\gamma_1-1} z_2^{\gamma_2}z_3^{\gamma_3+1}
\overline{z_1}^{\gamma_1} \overline{z_2}^{\gamma_2}\overline{z_3}^{\gamma_3}
\\

&+
\overline{a_{12}}c_R\overline{c_R}\gamma_1
{z_1}^{\gamma_1} {z_2}^{\gamma_2}{z_3}^{\gamma_3}
\overline{z_1}^{\gamma_1-1} \overline{z_2}^{\gamma_2+1}\overline{z_3}^{\gamma_3}
&+
\overline{a_{13}}c_R\overline{c_R}\gamma_1
{z_1}^{\gamma_1} {z_2}^{\gamma_2}{z_3}^{\gamma_3}
\overline{z_1}^{\gamma_1-1} \overline{z_2}^{\gamma_2}\overline{z_3}^{\gamma_3+1}.
\end{array}
\end{equation}
Without loss of generality,  it is enough to  study the term 
 $a_{12}c_P\overline{c_Q}\alpha_1
z_1^{\alpha_1-1} z_2^{\alpha_2+1}z_3^{\alpha_3}
\overline{z_1}^{\beta_1} \overline{z_2}^{\beta_2}\overline{z_3}^{\beta_3}.$ 
Using Proposition\eqref{LOG}, we conclude that the only possibility  for this term to  match with a term  containing $a_{13}$  is 
\begin{equation}
\begin{array}{ccccc}
\gamma_1&=&\alpha_1-1&=&\beta_1+1\\
\gamma_2&=&\alpha_2+1&=&\beta_2\\
\gamma_3&=&\alpha_3&=&\beta_3-1,
\end{array}
\end{equation}
with  $\gamma_1=0.$ This is a contradiction.

 \noindent Similarly, one can  show  that  terms containing  $a_{12}$ can  not interfer with terms containing  $a_{21}$  or  $a_{32}.$  We leave the details  to the reader.
 
 \noindent Therefore,  the     terms containing  $a_{12}$ can only interfere with  terms  containing  $a_{23}$ and $a_{31}.$
 
 \noindent Suppose that 
 $
a_{23}\neq0$ and that a term containing $a_{12}$ matches with a term $a_{23}.$  We have to study the   terms containing  $a_{12}$  and  $a_{23}$ given by

\begin{equation}
\begin{array}{lll}
&{\color{white}{+}}
a_{12}c_P\overline{c_Q}\alpha_1
z_1^{\alpha_1-1} z_2^{\alpha_2+1}z_3^{\alpha_3}
\overline{z_1}^{\beta_1} \overline{z_2}^{\beta_2}\overline{z_3}^{\beta_3}
&+
a_{23}c_P\overline{c_Q}\alpha_2
z_1^{\alpha_1} z_2^{\alpha_2-1}z_3^{\alpha_3+1}
\overline{z_1}^{\beta_1} \overline{z_2}^{\beta_2}\overline{z_3}^{\beta_3}
\\
&+
\overline{a_{12}}c_Q\overline{c_P}\alpha_1
{z_1}^{\beta_1} {z_2}^{\beta_2}{z_3}^{\beta_3}
\overline{z_1}^{\alpha_1-1} \overline{z_2}^{\alpha_2+1}\overline{z_3}^{\alpha_3}
&+
\overline{a_{23}}c_Q\overline{c_P}\alpha_2
{z_1}^{\beta_1} {z_2}^{\beta_2}{z_3}^{\beta_3}
\overline{z_1}^{\alpha_1} \overline{z_2}^{\alpha_2-1}\overline{z_3}^{\alpha_3+1}
\\
\\
&+
a_{12}c_Q\overline{c_P}\beta_1
z_1^{\beta_1-1} z_2^{\beta_2+1}z_3^{\beta_3}
\overline{z_1}^{\alpha_1} \overline{z_2}^{\alpha_2}\overline{z_3}^{\alpha_3}
&+
a_{23}c_Q\overline{c_P}\beta_2
z_1^{\beta_1} z_2^{\beta_2-1}z_3^{\beta_3+1}
\overline{z_1}^{\alpha_1} \overline{z_2}^{\alpha_2}\overline{z_3}^{\alpha_3}
\\
&+
\overline{a_{12}}c_P\overline{c_Q}\beta_1
{z_1}^{\alpha_1} {z_2}^{\alpha_2}{z_3}^{\alpha_3}
\overline{z_1}^{\beta_1-1} \overline{z_2}^{\beta_2+1}\overline{z_3}^{\beta_3}
&+
\overline{a_{23}}c_P\overline{c_Q}\beta_2
{z_1}^{\alpha_1} {z_2}^{\alpha_2}{z_3}^{\alpha_3}
\overline{z_1}^{\beta_1} \overline{z_2}^{\beta_2-1}\overline{z_3}^{\beta_3+1}
\\
\\
&+
a_{12}c_R\overline{c_R}\gamma_1
z_1^{\gamma_1-1} z_2^{\gamma_2+1}z_3^{\gamma_3}
\overline{z_1}^{\gamma_1} \overline{z_2}^{\gamma_2}\overline{z_3}^{\gamma_3}
&+
a_{23}c_R\overline{c_R}\gamma_2
z_1^{\gamma_1} z_2^{\gamma_2-1}z_3^{\gamma_3+1}
\overline{z_1}^{\gamma_1} \overline{z_2}^{\gamma_2}\overline{z_3}^{\gamma_3}
\\
&+
\overline{a_{12}}c_R\overline{c_R}\gamma_1
{z_1}^{\gamma_1} {z_2}^{\gamma_2}{z_3}^{\gamma_3}
\overline{z_1}^{\gamma_1-1} \overline{z_2}^{\gamma_2+1}\overline{z_3}^{\gamma_3}
&+
\overline{a_{23}}c_R\overline{c_R}\gamma_2
{z_1}^{\gamma_1} {z_2}^{\gamma_2}{z_3}^{\gamma_3}
\overline{z_1}^{\gamma_1} \overline{z_2}^{\gamma_2-1}\overline{z_3}^{\gamma_3+1}.
\end{array}
\end{equation}
Using Proposition \eqref{LOG} again, one can check that this implies
\begin{equation}
\begin{array}{ccccc}
\gamma_1&=&\alpha_1-1&=&\beta_1\\
\gamma_2&=&\alpha_2+1&=&\beta_2+1\\
\gamma_3&=&\alpha_3&=&\beta_3-1,
\end{array}
\end{equation}
with  $\gamma_1=0$  and  $\alpha_2=0$.
Therefore, we obtain 
\begin{equation}
\begin{array}{ccccc}
\gamma_1=0&&\alpha_1=1&&\beta_1=0\\
\gamma_2=1&&\alpha_2=0&&\beta_2=0\\
\gamma_3=k&&\alpha_3=k&&\beta_3=k+1,
\end{array}
\end{equation}
and  $M_d$ is given by 
\begin{equation}
\Im w=
c_P\overline{c_Q}z_1z_3^k\overline{z_3}^{k+1}
+
c_Q\overline{c_P}z_3^{k+1}\overline{z_1}\hspace{.15mm}\overline{z_3}^k
+
c_R\overline{c_R}
z_2z_3^k
\overline{z_2}\hspace{.15mm}\overline{z_3}^k,
\end{equation}
with
$$
N_1=a_{12}z_2\frac{\partial}{\partial z_1}
+a_{23}z_3\frac{\partial}{\partial z_2}  \in \5g_{0}.
$$

The case $
a_{31}\neq0$  is similar and  its proof is left to the reader.

\end{proof}
\noindent We now return to the proof of Theorem \eqref{L33}. 

\noindent  We notice that $N_1$ given by  Lemma \eqref{TheLemma}  is nilpotent  in all cases.

\noindent Let  $N$  be given by \eqref{lilia} with  $N\ne 0.$  After a  possible permutation of the variables, we may assume that $a_{12} \ne 0.$ If  $N_2:=N-N_1\ne 0,$ we apply again     Lemma\eqref{TheLemma} to $N_2.$ We  obtain that $N_2$ is nilpotent, with  $M_d$  given, up to permutations of the variables,  by 

    \begin{equation}	\label{lilia1}
 	\Im w=
	2\Re\left(
	c_P\overline{c_Q}
	z_1z_2^k
	\overline{z_2}^{k+1}
	\right)+
	c_R\overline{c_R}
	z_2^kz_3
	\overline{z_2}^k\overline{z_3},
	\end{equation} 
	 
\begin{equation}\label{lilia2}
	\Im w=
	2\Re\left(c_P\overline{c_Q}
	z_1z_3^{k}
	\overline{z_2}\hspace{.15mm}\overline{z_3}^{k}
	\right)	+
	c_R\overline{c_R}
	z_3^{k+1}
	\overline{z_3}^{k+1},
	\end{equation}
or	
\begin{equation}\label{lilia3}	\Im w=
	2\Re\left(c_P\overline{c_Q}
	z_1z_2^{\alpha_2}z_3^{\alpha_3}\overline{z_2}^{\alpha_2+1}\overline{z_3}^{\alpha_3}
	\right)	
	+
	c_R\overline{c_R}
	z_2^{\gamma_2}z_3^{\gamma_3}\overline{z_2}^{\gamma_2}\overline{z_3}^{\gamma_3}.\end{equation}
\end{proof}

\br
It is not hard to see that $D \in \5g_{0}$ given by \eqref{lilia0} depends on $3$ parameters. Indeed, one can decompose $D$ into $D= D_1 +D_2,$ where  $D_1 \in \5g_{0}$ has a real diagonal and $D_2 \in \5g_{0}$ has a purely imaginary diagonal.  Using Proposition\eqref{LOG}, one obtains that $D_1$ depends on $1$  real parameter and $D_2$ on $2$ real  parameters.
\er
\br \label{bey}
Inspecting the proof of Lemma\eqref{TheLemma}, we see that $N$ given by \eqref{lilia} depends on $3$ real parameters  if $M_d$ is given by \eqref{lilia1}, $2$ real parameters if $M_d$ is given by \eqref{lilia2} and $1$ real parameter if $M_d$ is given by \eqref{lilia3}. Otherwise, $N=0.$ 
\er
We illustrate the above results   with the following examples:
\be Let $M_d$ be given by 
\begin{equation}	\Im w=z_1z_3\overline{z_3}^2
	+{z_3}^2\overline{z_1}\hspace{.3mm}\overline{z_3}
	+z_2z_3\overline{z_2}\hspace{.3mm}\overline{z_3},\end{equation}
	Then $N$ is given by 
	\begin{equation}
	N
	=
	\begin{pmatrix} 
	0	&a+ib&ic \\
	0		 & 0&-a+ib\\
	0		 & 0&0
	\end{pmatrix},
	\end{equation}
	where 
$	a,b,c\in\mathbb{R}.$ Here, $N$ is nilpotent.
\ee
\be  Let $M_d$ be given by 
	\begin{equation}\Im w=iz_1z_2
	\overline{z_1}\hspace{.3mm}\overline{z_3}
	-i{z_1}{z_3}\overline{z_1}\hspace{.3mm}\overline{z_2}
	+z_1^2\overline{z_1}^2,\end{equation}
	
	Then $N$ is given by 
	\begin{equation}
	N
	=
	\begin{pmatrix} 
	0	&0&0 \\
	0		 & 0&a\\
	0		 & b&0
	\end{pmatrix}= \begin{pmatrix} 
	0	&0&0 \\
	0		 & 0&a\\
	0		 & 0&0
	\end{pmatrix}+ \begin{pmatrix} 
	0	&0&0 \\
	0		 & 0&0\\
	0		 & b&0
	\end{pmatrix},
\end{equation}	
	with
$	a,b\in\mathbb{R}.$ Here, $N$ is not nilpotent, but is the sum of $2$ nilpotents rotations.
\ee

\section{ The structure of $\5g_{c}$ for  the  $PQR$ problem }
In this section, we study the grading component 
 $\5g_c,$ that is, the set of   those  rigid  vector fields of weight  strictly bigger than $0$ (and  strictly less than one \cite{KMZ}).
 
 \noindent We have the following proposition:
\bp\label{key}
Let  $M_d$ be given by \eqref {jm13} with $P,$ $Q$ and $R$ non necessarily monomials,
	 and holomorphically nondegenerate. Then 
	\begin{equation}
	0\leq\dim \boldsymbol{\mathfrak{g}}_{\mathfrak{c}}\leq3.
	\end{equation}
	\ep

\begin{proof}
	Let  $X\in\boldsymbol{\mathfrak{g}}_{\mathfrak{c}}$. We have
	\begin{eqnarray}
	0&=&X(P)\overline{Q}+X(Q)\overline{P}+X(R)\overline{R}
	\nonumber\\\nonumber
	&&+Q\overline{X(P)}+P\overline{X(Q)}+R\overline{X(R)}.
	\end{eqnarray}
	The case $\deg P=\deg Q$ yields to 
	\begin{equation}
	X(P)\overline{Q}+X(Q)\overline{P}+X(R)\overline{R}=0.
	\end{equation}
	Using  Lemma\eqref{LOG1}, we get a contradiction. 
	
\noindent The   case  	$\deg P<\deg Q$
	yields to

	\begin{equation}X(Q)=0, \ 
	X(P)\overline{Q}+X(R)\overline{R}+Q\overline{X(P)}+R\overline{X(R)}=0.
	\end{equation}
	which implies 
	\begin{equation}X(P)=0,\end{equation}
	or
	\begin{equation}
	X(P)\overline{Q}+R\overline{X(R)}=0,
	\end{equation}
	or
	$$
	X(P)\overline{Q}+Q\overline{X(P)}=0\;\;\;\text{et}\;\;\;X(R)=0.
	$$
	The first case implies $X(R)=0,$ which is a contradiction.
	The second case implies
	
	\begin{eqnarray}
	X(P)&=&\alpha R,\nonumber\\
	X(R)&=&\beta Q\nonumber,\;\;\text{and}\\
	X(Q)&=&0,\nonumber
	\end{eqnarray}
	with $\beta=-\overline{\alpha}.$

	The third case implies
	\begin{eqnarray}
	X(P)&=&\alpha Q,\nonumber\\
	X(R)&=& 0\nonumber,\;\;\text{and}\\
	X(Q)&=&0,\nonumber
	\end{eqnarray}
	with $\alpha=ic$, $c\in\mathbb{R}.$ 
Using Cramer's rule, we obtain the desired conclusion.
	
\end{proof}

The following example shows that $\dim \boldsymbol{\mathfrak{g}}_{\mathfrak{c}}=3$ is reached.

\be~\label{exempleAvril22}
	Let $M_d$ be given by 
	\begin{equation}
	\Im w=z_1z_2z_3^2
	\overline{z_1}^2\overline{z_2}\hspace{0.15mm}\overline{z_3}^3
	+
	{z_1}^2{z_2}{z_3}^3
	\overline{z_1}\hspace{.15mm}\overline{z_2}\hspace{.15mm}\overline{z_3}^2
	+	z_1z_2z_3^3
	\overline{z_1}\hspace{0.15mm}\overline{z_2}\hspace{0.15mm}\overline{z_3}^3.
	\end{equation}
$M_d$ is holomorphically non degenerate since the Jacobian 
\begin{equation}
J(P,Q,R)	= \begin{vmatrix}
	P_{z_1}&P_{z_2}&P_{z_3}\\
	Q_{z_1}&Q_{z_2}&Q_{z_3}\\
	R_{z_1}&R_{z_2}&R_{z_3}
	\end{vmatrix}
	=
	\begin{vmatrix}
		z_2z_3^2&z_1z_3^2&2z_1z_2z_3\\
2{z_1}{z_2}{z_3}^3&{z_1}^2{z_3}^3&3{z_1}^2{z_2}{z_3}^2	\\
	  z_2z_3^3&z_1z_3^3&3z_1z_2z_3^2
	\end{vmatrix}=-z_1^3z_2^2z_3^7
	\end{equation}
	is not vanishing.
	
\noindent Let $X$ be  given by 
$	X=\sum_{j=1}^{3}f_j(z)\frac{\partial}{\partial z_j}.
	$
Inspecting the proof of Theorem\eqref{key}, we see that there are  $2$ cases to study. Either
	
		\begin{eqnarray}
	X(P)&=&\alpha R,\nonumber\\
	X(R)&=&\beta Q\nonumber,\;\;\text{and}\\
	X(Q)&=&0.\nonumber
	\end{eqnarray}
	or
		\begin{eqnarray}
	X(P)&=&\alpha Q,\nonumber\\
	X(R)&=&0\nonumber,\;\;\text{and}\\
	X(Q)&=&0.\nonumber
	\end{eqnarray}
	The first case leads to the following system:
	 \begin{equation}
	\begin{pmatrix}
	P_{z_1}&P_{z_2}&P_{z_3}\\
	R_{z_1}&R_{z_2}&R_{z_3}\\
		Q_{z_1}&Q_{z_2}&Q_{z_3}
	\end{pmatrix}
	\begin{pmatrix}
f_1\\
f_2\\
f_3
\end{pmatrix}
	=
	\begin{pmatrix}
		\alpha R\\
		\beta Q\\
		0
	\end{pmatrix}
\end{equation}
Using Cramer's rule, we obtain
\begin{eqnarray}
f_1(z)
&=&\frac{
\begin{vmatrix}
	\alpha z_1z_2z_3^3&z_1z_3^2&2z_1z_2z_3\\
\beta {z_1}^2{z_2}{z_3}^3&z_1z_3^3&3z_1z_2z_3^2\\
0&{z_1}^2{z_3}^3&3{z_1}^2{z_2}{z_3}^2
\end{vmatrix}
}{z_1^3z_2^2z_3^7}
=-\beta {z_1}^2,
\nonumber
\end{eqnarray}
\begin{eqnarray}
f_2(z)
&=&\frac{
	\begin{vmatrix}
z_2z_3^2&\alpha z_1z_2z_3^3&2z_1z_2z_3\\
z_2z_3^3& \beta {z_1}^2{z_2}{z_3}^3  &3z_1z_2z_3^2\\
2{z_1}{z_2}{z_3}^3&0&3{z_1}^2{z_2}{z_3}^2
\end{vmatrix}
}{z_1^3z_2^2z_3^7}=
3\alpha
z_2z_3
-
\beta
{z_1}z_2,
\end{eqnarray}
\begin{eqnarray}
f_3(z)
&=&\frac{
	\begin{vmatrix}
z_2z_3^2&z_1z_3^2&\alpha z_1z_2z_3^3\\
z_2z_3^3&z_1z_3^3& \beta {z_1}^2{z_2}{z_3}^3  \\
2{z_1}{z_2}{z_3}^3&{z_1}^2{z_3}^3&0
\end{vmatrix}
}{z_1^3z_2^2z_3^7}=
	-\alpha z_3^2
+
\beta {z_1}{z_3},
\nonumber
\end{eqnarray}
 where  $\alpha=a+ib$ et $\beta=-\overline{\alpha}=-a+ib.$

The second case leads to the  following solutions:

\begin{eqnarray}
f_1(z)
=
0,
\nonumber
\end{eqnarray}
\begin{eqnarray}
f_2(z)
=
3\alpha
z_1z_2z_3,\;\;\;\;\text{and}
\nonumber
\end{eqnarray}
\begin{eqnarray}
f_3(z)
=
-\alpha z_1z_3^2,
\nonumber
\end{eqnarray}
 $\alpha=ic, \ c \in \Bbb R.$
Hence $\dim_{\mathbb{R}}\boldsymbol{\mathfrak{g}}_{\mathfrak{c}}=3.$
\ee
   
\br
It is interesting to notice  that $\dim_{\mathbb{R}}\boldsymbol{\mathfrak{g}}_{\mathfrak{c}}=0$ if $n=1,$  $\dim_{\mathbb{R}}\boldsymbol{\mathfrak{g}}_{\mathfrak{c}}=1$ if $n=2,$ while here we have  an example with $\dim_{\mathbb{R}}\boldsymbol{\mathfrak{g}}_{\mathfrak{c}}=3.$
\er
\section{ The structure of $\5g_{-\frac{1}{d}},$  $\5g_{1 -\frac{1}{d}}$  and $\5g_{1}$ for  the  monomial  $PQR$ problem }
We now study the  remaining grading components of weights $-\dfrac{1}{d},$  ${1 -\dfrac{1}{d}}$  and ${1}.$  Indeed,  by Theorem 1.1 of \cite{KMZ}, there is one more left, the grading component of weight $-1,$  $\5g_{-1},$   whose  real dimension is $1,$  and  generated by $\partial_w.$  


\noindent We start with the following lemma
\bl \label{at}
Let  $M_d$ be given by \eqref {jm13} with $P,$ $Q$ and $R$ monomials,
	 and holomorphically nondegenerate. If $\5g_{-\frac{1}{d}} \ne \{0\}$ then $M_d$ is, after a possible permutation of the variables, given by  
	 \begin{equation}\label{jim13}
  M_d = \{(z,w) \in \mathbb C^{3} \times  \mathbb C \ | \  \Im w= z_1\bar Q +  Q\bar {z_1} +  R\bar R  \},
\end{equation}
with $Q$ and $R$ independent of $z_1.$
\el

\begin{proof}
Let $X \in \5g_{-\frac{1}{d}}$ be given by 
\begin{equation} X=a\frac{\partial}{\partial z_1}
	+b\frac{\partial}{\partial z_2}
	+c\frac{\partial}{\partial z_3}
	+g(z)\frac{\partial}{\partial w}.\end{equation}
	We first assume that $g \ne 0.$ It implies that $P,$ $Q$ or $R$ is linear. Without loss of generality, we may assume that  $P$ is    non zero linear since $R$ can not be linear  ($M_d$ is assumed to be Levi degenerate).  By permuting the variables, we may assume that $P= z_1,$  and hence $X=a\frac{\partial}{\partial z_1}.$ This implies
\begin{equation}\label{tel}
\Re (	X(Q)\bar {z_1} + X(R)\bar {R})=0.
\end{equation}
which is possible if $X(Q)=0,$ and $X(R)=0.$

\noindent We now assume that $g=0.$  By a holomorphic change of coordinates, we may assume that $X=\frac{\partial}{\partial z_1}.$
We obtain 
\begin{equation}
\Re (	X(P)\bar Q +  X(Q)\bar P +  X(R) \bar R)=0
\end{equation}
Considering the biggest power in $z_1$ in $P,$ $Q,$ $R, $ it is not hard to conclude that this case can not happen.
\end{proof}
	
We have the following proposition:

\bp \label{bientotfini}
	
	\nopagebreak[4]\noindent
	Let  $M_d$ be given by \eqref {jm13} with $P,$ $Q$ and $R$ monomials,
	 and holomorphically nondegenerate. Then 
	$$
	\dim\boldsymbol{\mathfrak{g}}_{-\frac{1}{d}}>0\;\;\;
	\Leftrightarrow\;\;\;
	\dim\boldsymbol{\mathfrak{g}}_{{1-\frac{1}{d}}}>0.
	$$
\ep

\bpf
Using Theorem 1.2 in  \cite{KM22}, we only have to show that there exists a holomorphic vector field $Y$
satisfying 
	\begin{equation}
Y (R) = QR, \ \ 
Y (Q) = Q^2.
\end{equation}
Since  $M_d$ is holomorphically nondegenerate and monomial, we obtain $Y$ explicitely, using Cramer's rule.
\epf

\bl\label{sara1}
Let  $M_d$ be given by \eqref {jm13} with $P,$ $Q$ and $R$ monomials,
	 and holomorphically nondegenerate. Then 	
	
\begin{equation}
	\dim_{\mathbb{R}}\boldsymbol{\mathfrak{g}}_{1}=1.
	\end{equation}
\el
\bpf
 By  Theorem 1.1 in \cite{KMZ}, we obtain that $\dim_{\mathbb{R}}\boldsymbol{\mathfrak{g}}_{1}\ne 0$ if and only if there exists a vector field $Y$ with $$Y(P\bar Q +Q\bar P +R\bar R)=P\bar Q +Q \bar P +R \bar R.$$  Using Lemma \eqref{great1}, we obtain that $Y$ is real diagonal. We   then solve  a Cramer system, and   using Lemma \eqref{LOG}, get explicitely $Y.$	 
\epf

\section{ Statements and Proofs of the main results}

\bt \label{MainResult}
	Let  $M_d$ be given by \eqref {jm13} with $P,$ $Q$ and $R$ monomials,
	 and holomorphically nondegenerate. Assume that $\dim \boldsymbol{\mathfrak{g}}_c >0.$ Then 
	 
\begin{equation}
	7\leq\dim \boldsymbol{\mathfrak{g}}\leq13.
	\end{equation}
	Furthermore,  there exists a hypersurface $M_d$ for which $ \dim \boldsymbol{\mathfrak{g}}=13.$
\et


\bpf
Using the characterization of the existence of a nilpotent rotation given by \eqref{lilia1}, \eqref{lilia2} and \eqref{lilia3}, and inspecting the proof of Theorem \eqref{key}, we conclude that under the assumption $\dim \boldsymbol{\mathfrak{g}}_c >0,$ there are no nilpotent rotations. Using  Lemma \eqref{at},  Proposition \eqref{bientotfini}, Lemma  \eqref{sara1} and Theorem \eqref{key}, we conclude that $7\leq\dim \boldsymbol{\mathfrak{g}}\leq13.$

\noindent  It is not hard to see that  the hypersurface given by $$\Im w=  z_1 \overline {{z_2}}^l+  z_2^l \overline {z_1} +z_2^{\frac{\ell-1}{2}}z_3
	\overline{z_2}^{\frac{\ell-1}{2}}\overline{z_3}$$  admits $\dim \boldsymbol{\mathfrak{g}}_c =3.$
		
Indeed, the following vectors fields  are in $\boldsymbol{\mathfrak{g}}_c:$
\begin{enumerate}	
		\item $X_1=aiz_2^\ell\frac{\partial}{\partial z_1}$, $a\in\mathbb{R}^,$
		
\item $X_2= - (a+ib)z_2^{\frac{\ell-1}{2}}z_3\frac{\partial}{\partial z_1}
		+ (a-ib)z_2^{\frac{\ell+1}{2}}\frac{\partial}{\partial z_3},$ $a, b\in\mathbb{R}.$
		\end{enumerate}
		Hence, Using  Lemma \eqref{at} and   Proposition \eqref{bientotfini}, we conclude that $\dim \boldsymbol{\mathfrak{g}}=13.$
 \epf

\bt \label{MainResult2}
	Let  $M_d$ be given by \eqref {jm13} with $P,$ $Q$ and $R$ monomials,
	 and holomorphically nondegenerate. Assume that $\dim \boldsymbol{\mathfrak{g}}_c =0.$ Then 
	 
\begin{equation}
	6\leq\dim \boldsymbol{\mathfrak{g}}\leq 9.
	\end{equation}
	Furthermore,  there exists a hypersurface $M_d$ for which $ \dim \boldsymbol{\mathfrak{g}}=9.$ 
\et	

\bpf  If $\dim \boldsymbol{\mathfrak{g}}_c =0,$ then
 by  Lemma \eqref{at}, $\dim \boldsymbol{\mathfrak{g}}_{-\frac{1}{d}}  =0.$ 
 Therefore,  we have  $	6\leq\dim \boldsymbol{\mathfrak{g}}\leq 9.$  Using Remark \eqref{bey}, we conclude that $\dim \boldsymbol{\mathfrak{g}} =9$ for the hypersurface given by \eqref{lilia1}.

  \epf


\begin{thebibliography}{BER96b}

\itemsep=2pt

%
\bibitem{Be}
V.K.Beloshapka, V.V. Ezhov, G. Schmalz, {\it Holomorphic classification of four-dimensional surfaces in $\Bbb C^3$,  }  Izv. Ross. Akad. Nauk Ser. Mat.  \textbf{72}
(2008),  3--18.  
\bibitem{BG}
Bloom, T.,  Graham, I., {\it On "type"  conditions for
generic real submanifolds of $  C\sp{n}$},  Invent. Math.  \textbf{40}
(1977),  217--243.

\bibitem{BFG} Beals, M., Fefferman, C., Graham R.,
{\it Strictly pseudoconvex domains in $\mathbb C^n$},
Bull.\ Amer.\ Math.\ Soc.\ (N.S.)
 \textbf{8} (1983), 125--322.
 

 

 \bibitem{C2}  E.Cartan : \textit{Sur la g\'eom\'etrie pseudo-conforme
des hypersurfaces de deux variables complexes, II}, Ann.Scoula
Norm. Sup. Pisa   \textbf{1} (1932), p. 333--354.



\bibitem{CM} Chern, S.\ S.\, Moser, J., \textit{Real hypersurfaces in
complex manifolds},  Acta Math.\  \textbf{133} (1974),  219--271.










\bibitem{IK}
Isaev, A. V., Kruglikov, B., {\it On the symmetry algebras of 5-dimensional CR-manifolds}, Adv. Math. {\bf 322} (2017), 530--564.


\bibitem{KK} Kim, S. Y.,  Kol\'a\v r, M.,  \textit{Infinitesimal symmetries of weakly pseudoconvex manifolds
}, Math. Z.,  \textbf{300} (2022), no. 3, 2451--2466.

\bibitem{K} Kohn, J.\ J., \textit{Boundary behaviour of
$\bar \partial$ on weakly pseudoconvex manifolds of dimension
two},
 J.\ Differential  Geom.\  \textbf{6} (1972),  523--542.





\bibitem{Ko1a} Kol\'a\v r, M., {\it Normal forms for hypersurfaces of finite type in $ \mathbb C^2$}, Math. Res. Lett., \textbf{12} (2005),  523--542.


\bibitem{FM17}
Kol\'a\v r, M., Meylan, F., Infinitesimal CR automorphisms   for a class of  polynomial models
Arch. Math. (Brno) {\bf 53} no. 5 (2017): 255-265
\bibitem{FM}
Kol\'a\v r, M., Meylan, F., {\it Chern-Moser operators and weighted jet determination problems},
Geometric analysis of several complex variables and related topics, 75--88, Contemp. Math. 550, 2011.
\bibitem{KoMe}
Kol\'a\v r, M., Meylan, F. \textit{ Infinitesimal CR automorphisms of hypersurfaces of finite type in C2.},  Arch. Math. (Brno) {\bf 47} no. 5 (2011): 367-375
\bibitem{KMZ} M. Kol\'a\v r, F.  Meylan, D. Zaitsev, {\it Chern-Moser operators and polynomial models in CR geometry},  Adv. Math.  \textbf{263} (2014), 321-356.
\bibitem{KM22}Kol\'a\v r, M., Meylan, F.,\textit{ Nonlinearizable CR Automorphisms for Polynomial Models in $\Bbb C^N.$}, J. Geom. Anal. {\bf 33} no. 3(2023): 106-






\bibitem{KL}
Kruzhilin, N. G., Loboda, A. V.,
 {\it Linearization of local automorphisms of pseudoconvex surfaces},
 Dokl. Akad. Nauk SSSR,
\textbf{271} (1983), 280--282.

\bibitem{Kr20} Kruglikov B., {\it  
Blow-ups and infinitesimal automorphisms of CR-manifolds}
Mathematische Zeitschrift  \textbf{296} (2020), 1701-1724.




\bibitem{Po} Poincar\'e, H.,  \textit{Les fonctions analytiques de
deux variables et la repr\'esentation conforme, } Rend. Circ. Mat.
Palermo \textbf{23} (1907), 185--220.


\bibitem{S} Stanton, N.,
{\it Infinitesimal CR automorphisms of real hypersurfaces}, Amer. J. Math. {\bf 118} (1996),  209--233.
\bibitem{Tanaka} Tanaka N. I.,
	\emph{On the pseudo-conformal geometry of hypersurfaces of the space of $n$
		complex variables}.
	J. Math. Soc. Japan, \textbf{14} (1962), pp.397-429.
\bibitem{Tr} Treves, F., {\it A treasure trove of geometry and analysis: the hyperquadric}, Notices Amer. Math. Soc. \textbf{47} (2000),  

\bibitem{V} Vitushkin, A.G.,  \textit{Real analytic
hypersurfaces in complex manifolds}, Russ. Math. Surv. \textbf{40}
(1985),  1--35.

\bibitem{W} Webster, S.M., \textit{On
the Moser normal form at a non-umbilic point}, Math. Ann.
 \textbf{233} (1978),  97--102.


\bibitem{Y} Yang, P., {\it Automorphism of tube domains}, Amer. J. Math. {\bf 104} (1982), 1005--1024. 

\end{thebibliography}
\end{document}